\documentclass[a4paper,11pt]{article}
\usepackage[a4paper, margin=2.25cm]{geometry}
\usepackage[bitstream-charter]{mathdesign}
\usepackage{changes}
\usepackage{amsmath}
\usepackage{url}
\usepackage{enumitem}
\usepackage{booktabs}
\usepackage{listings}
\usepackage{xcolor}
\usepackage{authblk}
\usepackage{marvosym}
\usepackage{framed}
\usepackage{graphicx}
\usepackage{verbatimbox}
\usepackage{multirow}
\usepackage{gensymb}
\usepackage[colorlinks=true, allcolors=blue, breaklinks=true]{hyperref}
\usepackage{tikz}
\usetikzlibrary{arrows.meta}
\usepackage{cite}
\usepackage{amsthm}
\usepackage{kotex}
\usepackage{textcomp}
\DeclareMathOperator*{\Res}{Res}
\theoremstyle{definition}

\newtheorem{example}{Example}

\makeatletter
\newcommand{\tpmod}[1]{\!\!\!{\@displayfalse\pmod{#1}}}
\makeatother
\renewcommand{\arraystretch}{1.5}
\newcommand{\atan}{\operatorname{atan}}
\title{\vspace*{-1cm}\bfseries \Large Calculus and digital natives in rendezvous: \emph{wxMaxima} impact}
\author{Natanael Karjanto\thanks{\Letter: \texttt{natanael@skku.edu} \; \href{https://orcid.org/0000-0002-6859-447X}{\includegraphics[scale=0.08]{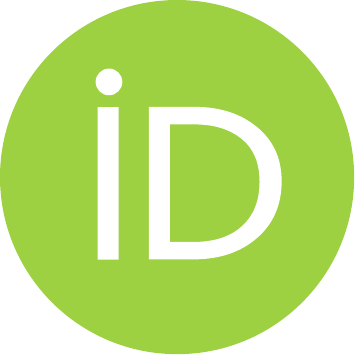}}}}
\affil{Department of Mathematics, University College, Natural Science Campus, Sungkyunkwan University, Suwon~16419, Republic of Korea}
\date{\vspace*{-1cm} \scriptsize Last updated \today}

\begin{document}
\maketitle 
\vspace*{-1cm}
\begin{abstract}
{\noindent
This article covers how a computer algebra system (CAS) \emph{wxMaxima} can be explored for teaching and learning Single-Variable and Multivariable  Calculus for Korean digital natives. We present several examples where \emph{wxMaxima} can handle Calculus problems easily, not straightforwardly but still successfully solved with some human intervention, and unsuccessfully. By soliciting qualitative feedback on students' experience in exploiting the CAS, we gathered a mixed reaction. Although some students commented positively, the majority seemed to be resistant in embracing a new technological tool. \vspace*{0.2cm} }
\end{abstract}

\section{Introduction}

Each year, college freshmen in South Korea who plan to major in Science and Engineering are required to enroll in two Calculus courses: Single-Variable and Multivariable Calculus, hereby SVC and MVC, respectively. Growing up surrounded by technology, they belong to a group called ``digital natives''. Meanwhile, teaching with technology is not only essential but also has become a necessity during the Fourth Industrial Revolution. With so much software to choose from, we selected a free computer algebra system (CAS) \emph{wxMaxima}. In this article, we will discuss some features of this software, particularly in connection to symbolic computation for teaching SVC and MVC. After introducing the CAS into the classrooms, we discuss the students' feedback and reflect on their opinion to improve the pedagogical approach. In this introduction, we will briefly cover the main components in this study: Calculus, digital natives, and CAS \emph{wxMaxima}. Prior to embedding and implementing \emph{wxMaxima} into our classrooms, we hypothesize that Korean digital natives would embrace the CAS more positively than migrant, non-digital natives. It turned out that our premise is not entirely correct.

\subsection{What is Calculus?}

Calculus is the branch of mathematics that studies how phenomena change, how to measure that change, and we utilize those measurements in our life. It administers the foundation for modeling systems where changes are present. The subject was invented because earlier scientists and mathematicians had a great interest in physical sciences. In the past, Calculus was taught after secondary school level, but nowadays, Calculus is also taught and introduced in high school. By the 1930s, Calculus became an important part of high school mathematics in the US~\cite{swenson}. By the 1960s, the idea of not teaching Calculus in high school would have been unthinkable~\cite{ferguson}.

Although Calculus might be the peak of high school mathematics and many students are rushing into taking it, many educators and mathematicians are concerned that there are some drawbacks of such movement. They argued that many students who have taken Calculus in high school and enrolled in college-level Calculus courses tend to view the subject as another repetition and thus spend less time in understanding the concepts. Eventually, some struggle and might lose interest in mathematics entirely~\cite{bressoud,krantz}. What should be an exciting course, Calculus was branded as a barrier or filter course, particularly for STEM (Science, Technology, Engineering, and Mathematics) majors, instead of a pump~\cite{moore,subramaniam2008,hieb,bego,steen}.

In many universities across North America, the sequence of Calculus courses is often offered as a part of the undergraduate curriculum, particularly during freshman year. SVC is usually split into two-semester, one-year-long courses. The so-called ``Calculus~1'' usually covers Differential Calculus while ``Calculus~2'' comprises Integral Calculus and most often encompasses an Introduction to Differential Equations and Sequence and Series. What is known as ``Calculus~3'' commonly deals with Multivariable or Vector Calculus, and is usually offered during sophomore year.

In South Korea, on the other hand, both SVC and MVC are compressed into the freshmen year. During their first semester in college, i.e., Spring (Northern Hemisphere), the freshmen take the so-called ``Calculus~1'', which comprises both Differential and Integral Calculus, the combined course materials from what is commonly offered in North American universities as two separate courses. After the Summer break, the second semester of freshmen (Fall, Northern Hemisphere), freshmen enroll in ``Calculus~2'', which is equivalent to ``Calculus~3'' (MVC) offered by North American universities.

\subsection{Who are digital natives?}

The majority of our students, if not all, are ``digital natives'' or ``net generation''. A term coined by Marc Prensky, ``digital native'' is used to describe the generation of people who grew up in the era of ubiquitous technology, which includes computers and the Internet~\cite{prensky2001a,prensky2001b}. Helsper and Eynon argued that we often erroneously presume our ability as educators to engage our students with technology due to the existence of a generational gap. It turns out that this is only one factor; other factors could be more important than generational differences, including gender, experience, educational level, and the intensity of use~\cite{helsper}. 

Indeed, while technology is embedded in young people's lives, their use and skills are not uniform~\cite{bennet}, cf.~\cite{bennet2010beyond}. A similar finding among Australian university freshmen confirmed this lack of homogeneity concerning technology~\cite{kennedy}. Another study showed that in addition to a limited range of established technologies being utilized by young people, there was no evidence that they adopt radically different learning styles~\cite{margaryan}.

Although South Korea is one of the world's most technologically advanced and digitally connected countries, several case studies do indeed support research findings from other countries. For example, Koh and Shin observed a limited degree and homogeneous pattern in media use among Korean youths~\cite{koh2018}. On the other hand, Chung discovered the heterogeneous pattern in multi-tasking behavior and most Korean college students exhibit an affection-type of dependency on digital technology~\cite{chung2019}.

A literature review on digital natives presented in this article is far from exhaustive. For more extensive coverage on the topic, the following three references will be helpful. Kivunja provided a literature review to shed some light on theoretical perspectives of how digital natives learn and how we can use that knowledge in facilitating their learning~\cite{kivunja}. Palfrey and Gasser offered a sociological portrait of the first generation of digital natives who can seem both extraordinarily sophisticated and strangely narrow~\cite{palfrey}. Dingly and Seychell extensively covered the topics of digital natives, investigated the paradigm shift between the different generations of digital natives, and analyzed the future trends of technology~\cite{dingli}.

\subsection{Teaching with technology}

At the beginning of each new academic year, we as teachers and instructors will teach younger students than the previous cohorts. Since our digital natives tend to attach to technology even at a younger age than earlier generations, teaching with technological tools is no longer an option but has become a necessity. Granted, the current prolonged COVID-19 pandemic has only accelerated this transition process. In particular, teaching and learning mathematics using CAS has intensified during the past two decades. This subsection provides a brief overview of some studies where other software have been embedded in the mathematics classroom. The following subsection will cover \emph{wxMaxima} in particular and why we opt for it.

\emph{GeoGebra} is a dynamic mathematics software for all levels of education. As the name suggests, it not only brings together geometry and algebra but also calculus, statistics, and graphing in a single easy-to-use package. Hohenwarter et al. presented applications of \emph{GeoGebra} for teaching Calculus at both secondary and college levels~\cite{hohenwarter}. The findings from Saha et al. suggest that using \emph{GeoGebra} enhanced students' academic performance in understanding coordinate geometry~\cite{saha}. Other studies offer insights on how \emph{GeoGebra} generated not only enjoyment and fun in learning mathematics but also in concretizing abstract concepts~\cite{celen}. Even though utilizing \emph{GeoGebra} in the mathematics classroom still poses challenges and limitations~\cite{wassie}, the software has been found particularly useful in enhancing flipped classroom pedagogy~\cite{weinhandl}.

Another well-known CAS is \emph{Maple}, developed by a Canadian-based software company \emph{Maplesoft}. The selection of the name should not come as a surprise since the maple leaf is the most widely recognized national symbol of Canada. \emph{Maple} has been incorporated in project-based learning for Calculus class~\cite{wu}. The software's implementation in the Calculus classroom has been arguably positive~\cite{ningsih,lestiana,hamid,purnomo}. Even though many teachers and college instructors struggle in delivering online teaching during the COVID-19 pandemic, an initiative from Italy bears a fruitful resolution. By combining \emph{Maple} with \emph{Moodle} digital learning environment, teachers employ \emph{Maple} in problem-solving activities, designing student worksheets, increasing students' participation, and providing immediate feedback~\cite{fissore}.

A rigorous competitor of \emph{Maple} and \emph{wxMaxima} is \emph{Wolfram Mathematica}, or simply \emph{Mathematica}. The software was formulated by Stephen Wolfram and is currently being developed by Wolfram Research in Champaign, Illinois. Although it has been criticized for being a closed source software, it possesses the capabilities of high-performance computing. \emph{Mathematica} is integrated with \emph{WolframAlpha}, a computational knowledge engine developed by the same company. Dimiceli et al. described both benefits and drawbacks in teaching Calculus using \emph{WolframAlpha}~\cite{dimiceli}. Barba-Guaman et al. utilized \emph{Mathematica} to improve reading comprehension in mathematics~\cite{barba}. Beyond Calculus, \emph{Mathematica} enhances students' creativity and academic performance in Linear Algebra~\cite{rahmawati}, Discrete Mathematics~\cite{ivanov}, Classical Mechanics~\cite{romano}, and Economics~\cite{rihova}, among others. Table~\ref{CAS-summary} summarizes the literature study on teaching with technology.
\setlength{\tabcolsep}{0.25cm}
\renewcommand{\arraystretch}{1.0}
\begin{table}[h]
\begin{center}
\begin{tabular}{@{}llll@{}}
\toprule
Software 			& Creator (Started)      & Literature 								& Topic \\ \hline
\multirow{4}{*}{\emph{GeoGebra}}	& 		 & Hohenwarter et al.~\cite{hohenwarter} 	& Calculus \\
					& Markus Hohenwarter	 & Saha et al.~\cite{saha}		 			& Coordinate geometry \\										
					& (2001)				 & Celen~\cite{celen}		 				& Line and angle \\
					&						 & Wassie and Zergaw~\cite{wassie}		 	& Precalculus \\ \hline
\multirow{6}{*}{\emph{Maple}}	&			 & Wu and Li~\cite{wu} 						& Calculus \\
                    & 		                 & Ningsih and Paradesa~\cite{ningsih}		& Freshmen college math \\
                    & University of Waterloo & Lestiana and Oktanivani~\cite{lestiana}  & Integral calculus \\
                    & (1980)				 & Hamid et al.~\cite{hamid}			    & Integral calculus \\
                    &						 & Purnomo et al.~\cite{purnomo}			& Multivariable calculus \\
                    &						 & Fissore et al.~\cite{fissore}			& Contextualized problems \\ \hline	
\multirow{6}{*}{\emph{Mathematica}}	&  		 & Dimiceli et al.~\cite{dimiceli} 			& Calculus \\
                    & 						 & Barba-Guaman et al.~\cite{barba}			& Computer science \\
					& Wolfram Research		 & Rahmawati et al.~\cite{rahmawati}		& Linear algebra \\
                    & (1986)				 & Ivanov et al.~\cite{ivanov}				& Discrete mathematics \\
                    &						 & Romano and Marasco~\cite{romano}			& Classical mechanics \\
                    &						 & Říhová et al.~\cite{rihova}				& Economics \\ \hline	
\multirow{8}{*}{\emph{Maxima}}	&			 & Hannan~\cite{hannan}						& Calculus \\
				 	&       		         & Díaz et al.~\cite{diaz}					& Linear Algebra \\
				 	&						 & Ayub et al.~\cite{ayub}					& Secondary mathematics \\
				 	& Bill Schelter et al.	 & Dehl~\cite{dehl}							& Vector calculus \\
				 	& (1976)				 & Timberlake and Mixon~\cite{timberlake}	& Classical mechanics \\
				 	&						 & Senese~\cite{senese}						& Chemistry \\
				 	&						 & Žáková~\cite{zakova}						& Engineering \\
				 	&						 & Fedriani and Moyano~\cite{fedriani}		& Miscellaneous \\		
\bottomrule
\end{tabular}
\end{center}
\caption{A list of selected literature studies where CAS has been assimilated to teaching and learning, not only in mathematics and calculus but also in other subjects.}	\label{CAS-summary}
\end{table}

\subsection{Why do we opt for \emph{wxMaxima}?}

\emph{wxMaxima} is one of the graphical user interfaces (GUI) for \emph{Maxima}. \emph{Maxima} is one of CAS that specializes in symbolic computation. A CAS is software that can solve computational problems by rearranging formulas and providing analytical expressions instead of spitting out numerical values. Both systems are free of charge and released under the terms of the GNU General Public License (GPL). Under this authorization, everyone has the right and freedom to modify and distribute the software as long as its license with them remains unmodified.

Some of the free-of-charge competitors of \emph{wxMaxima} are \emph{Axiom} and \emph{SageMath}. The former lacks a GUI and the latter integrates many CAS packages into a common GUI using a syntax resembling \emph{Python}. Although \emph{SageMath} provides a cell server, the results are often slow to appear. Several commercial programs are more famous than \emph{wxMaxima} and we need to admit that \emph{wxMaxima} might not be able to compete with them. The so-called \emph{4M} CAS are developed by companies that employ programmer and developer experts: \emph{Magma}, \emph{Maple}, \emph{Mathematica}, and \emph{Matlab}. The latter initially specialized in numerical computations but the symbolic tools have been added later.

For mathematics teachers and college instructors who seek a free alternative for CAS, adopting \emph{wxMaxima} might be worth a try. In addition to serving as a calculator and symbolic manipulator, it can sketch two- and three-dimensional objects with high-quality figures. Since the CAS is lightweight and straightforward, performing simple computations can be completed swiftly. However, a word of precaution should be taken into account, as we will see in Section~\ref{symbolic}. Some problems need modifications for \emph{wxMaxima} to be able to solve, while other problems cannot be solved at all. Indeed, \emph{wxMaxima} possesses many limitations, particularly when it comes to symbolic integration.

The literature offers some glimpse on the CAS. Some promising attempts have been shown for teaching Calculus and Linear Algebra using \emph{wxMaxima}~\cite{hannan,diaz,ayub,dehl}. In Physics, Chemistry, Engineering, and even Business, the software is also starting gaining popularity~\cite{timberlake,senese,zakova,fedriani}. Evaluating the impact of symbolic computation in education has also been addressed~\cite{li}. In this article, our focus is on the teaching and learning of Calculus. It extends our previous discussion and complements the mathematical facet of \emph{wxMaxima} with educational aspects~\cite{karjanto17,karjanto21}.
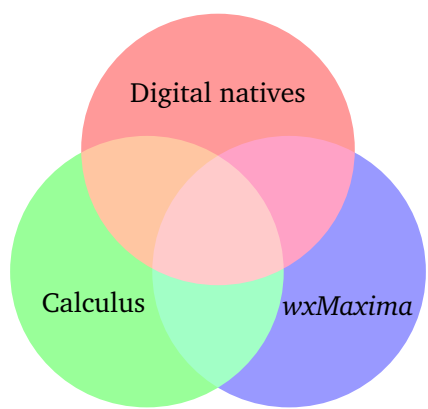
\begin{figure}[h]
\begin{center}
\begin{tikzpicture}[scale=0.9]
\begin{scope}[blend group = soft light]
\fill[red!40!white]   ( 90:1.2) circle (2);
\fill[green!40!white] (210:1.2) circle (2);
\fill[blue!40!white]  (330:1.2) circle (2);
\end{scope}
\node at ( 90:2)    {Digital natives};
\node at ( 210:2.1)   {Calculus};
\node at ( 330:2.2)   {\emph{wxMaxima}};
\end{tikzpicture}
\end{center}
\caption{A theoretical framework of teaching Calculus to digital natives using CAS \emph{wxMaxima}.}	\label{fig1-framework}
\end{figure}
\vspace*{-0.75cm}

\subsection{Theoretical framework and research question}

From the preceding subsections, we have collected three components so far. The object: Calculus courses (Single Variable and Multivariable), the subject study (Korean digital natives), and the pedagogical tool (CAS \emph{wxMaxima}). Due to the advancement of technology and the majority of our students are digital natives, then teaching Calculus using technology to the net generation is no longer an option; it becomes a necessity. By combining these three components, we obtain a Venn-diagram-like relationship. In turn, their intersection forms our theoretical framework for this study, as shown in Figure~\ref{fig1-framework}.

We are interested in investigating the following research questions: \vspace*{-0.5cm}
\begin{enumerate}[leftmargin=1.5em]
\setlength\itemsep{0pt}	
\item What type of symbolic computation that \emph{wxMaxima} can perform easily, with some manipulation, and none? 
\item After adopting and implementing \emph{wxMaxima} into Calculus teaching and learning, what kind of feedback do we receive as instructors? And what can we learn from these students' feedback?
\end{enumerate}

This paper is organized as follows. After this introduction, Section~\ref{symbolic} discusses symbolic computation using \emph{wxMaxima}. Through several examples, we will observe what \emph{wxMaxima} can solve easily and quickly, a type of problem that needs human intervention for \emph{wxMaxima} to be able to solve, and another category of problems where \emph{wxMaxima} is clueless. Section~\ref{methodology} continues with the educational aspect of \emph{wxMaxima}. It covers research methodology in how we collect data of students' feedback after embedding the CAS into the Calculus classrooms. Some findings are also presented. Section~\ref{conclusion} concludes this study.

\section{Symbolic computation using \emph{wxMaxima}}		\label{symbolic}

For simple problems, symbolic computation using \emph{wxMaxima} can be performed within seconds. In this section, we consider examples when \emph{wxMaxima} is not only reliable but also struggles in providing immediate outputs. The following examples provide a brief overview of the kinds of computations \emph{wxMaxima} can and cannot do.

\begin{example}[Easily solved]
Taylor and Maclaurin series of a function
\begin{equation*}
y = f(x) = 2x e^{-x^2} \cos(3x).
\end{equation*}
\end{example}
Finding Taylor and Maclaurin series are exceptionally fast. The \emph{wxMaxima} command\\ ``\verb|taylor(f(x),x,%pi,2);|'' gives the first three nonzero terms of the Taylor series expansion for the function~$f$:
\begin{equation*}
e^{-\pi^2} \left[-2\pi + (4\pi^2 - 2)(x - \pi) - (4\pi^3 - 7 \pi) (x - \pi)^2 + \dots \right].
\end{equation*}
By replacing $\pi$ with $0$, the syntax ``\verb|taylor(2*x*%e^(-x^2)*cos(3*x),x,0,8);|'' gives the first four nonzero terms of the Maclaurin series expansion:
\begin{equation*}
2x - 11 x^3 + \frac{67}{4} x^5 - \frac{1633}{120} x^7 + \dots.
\end{equation*}
Admittedly, deriving these series by hand calculation can be tedious. Although every Calculus student should master it, \emph{wxMaxima} can be used to check the results or to compute Taylor series for extended functions that would be too arduous for hand computations. Easing the computational load allows the learners to focus more deeply on concepts and patterns.

Taylor and Maclaurin polynomials are employed to approximate functions using the partial sums of the corresponding Taylor and Maclaurin series, respectively. They extend the idea of linearization and are advantageous in understanding asymptotic behavior, the growth of functions, evaluating definite integrals, and solving differential equations. 

When it comes to applications, Taylor polynomials have an abundance of them. For example, an equation describing refraction at a spherical interface can be simplified by either a linear or quadratic approximation of the position angle variable. The former is known as the first-order optics or Gaussian optics and the latter is known as the third-order optics. This resulting optical theory has become the basic theoretical tool used to design lenses~\cite{hecht}.

Evaluating integrals of a rational function can often be easier when the function is decomposed and expressed in a combination of simpler rational functions. This is the idea of partial fraction decomposition.
\begin{example}[Easily solved with slightly different outputs]
Partial fraction decomposition of a rational function
\begin{align*}
P(x) = \frac{2x^2 + 21x - 36}{x^7 + 10 x^6 + 37 x^5 + 60 x^4 + 36x^3}.
\end{align*}
\end{example}
Using a simple command ``\verb|P:(2*x^2+21*x-36)/(x^7+10*x^6+37*x^5+60*x^4+36*x^3);|'', ``\verb|P1:partfrac(P,x)|'', we obtain the following instant result
\begin{equation*}
P_1(x) = \frac{26}{3 \left( x+3\right) }+\frac{3}{{{\left( x+3\right) }^{2}}}-\frac{6}{x+2}+\frac{35}{4 {{\left( x+2\right) }^{2}}}-\frac{8}{3 x}+\frac{9}{4 {{x}^{2}}}-\frac{1}{{{x}^{3}}}.
\end{equation*}
Observe that we can factor the denominator as $x^3 (x^2 + 5x + 6)^2$, which explains why we obtain the terms $x^m$, $(x + 2)^{n}$, and $(x + 3)^n$, for $m = 1,2,3$ and $n = 1, 2$ in the denominators of the partial fraction decomposition.
Integrating the original rational function (``\verb|integrate(P,x);|'') and its partial fraction decomposition (``\verb|integrate(P1,x);|'') using \emph{wxMaxima} produces slightly different outputs, given as follows respectively:
\begin{align*}
\int P(x) \, dx &= \frac{26 \log{\left( x+3\right) }}{3}-6 \log{\left( x+2\right) }-\frac{8 \log{(x)}}{3}-\frac{14 {{x}^{3}}+43 {{x}^{2}}+11 x-3}{{{x}^{4}}+5 {{x}^{3}}+6 {{x}^{2}}} \\
\int P_1(x) \, dx &= \frac{26 \log{\left( x+3\right) }}{3}-6 \log{\left( x+2\right) }-\frac{8 \log{(x)}}{3}-\frac{3}{x+3}-\frac{35}{4 \left( x+2\right) }-\frac{9}{4 x}+\frac{1}{2 {{x}^{2}}}.
\end{align*}
It should be noted that for rational functions that cannot be decomposed, \emph{wxMaxima} often fails to integrate such expressions.

\begin{example}[Require some manipulations]
Integral of a beta function
\begin{equation*}
\int \frac{1}{(x^2 - x^3)^{\frac{1}{3}}} \, dx \qquad \qquad \text{and} \qquad \qquad \int_{0}^{1} \frac{1}{(x^2 - x^3)^{\frac{1}{3}}} \, dx
\end{equation*}
\end{example}
\emph{wxMaxima} fails to evaluate both these indefinite and definite integrals. The commands \linebreak ``\verb|integrate(1/(x^2-x^3)^(1/3),x);|'' and ``\verb|integrate(1/(x^2-x^3)^(1/3),x,0,1);|'' \linebreak produce no output. However, by pulling out the factor $x^{2/3}$ from the integrand and rewriting it as a beta function, using the commands ``\verb|integrate(1/(x^(2/3)*(1-x)^(1/3)),x);|'' and \linebreak ``\verb|integrate(1/(x^(2/3)*(1-x)^(1/3)),x,0,1);|'', we obtain a result for both integrals:
\begin{align*}
\int \frac{1}{x^{\frac{2}{3}}(1 - x)^{\frac{1}{3}}} \, dx &= \ln \left[\frac{\sqrt[3]{x} + \sqrt[3]{1 - x}}{\sqrt[3]{x}} \right] - 
\frac{1}{2} \ln \left[\frac{x - \sqrt[3]{x^{2}}\sqrt[3]{1 - x} + \sqrt[3]{x} \sqrt[3]{(1 - x)^2}}{x} \right] \\ 
& \quad - \sqrt{3} \atan \left[\frac{2\sqrt[3]{1 - x} - \sqrt[3]{x}}{\sqrt{3} \sqrt[3]{x}} \right] \\
\int_0^1 \frac{1}{x^{\frac{2}{3}}(1 - x)^{\frac{1}{3}}} \, dx &= \frac{2\pi}{\sqrt{3}} = \text{beta} \left(\frac{1}{3}, \frac{2}{3} \right) = \frac{\pi}{\sin \frac{\pi}{3}}.
\end{align*}
To find the numerical value for the latter, the command ``\verb|float((2*%pi)/sqrt(3));|'' yields \linebreak $3.627598728468436$.

The beta function, which is closely related to the gamma function, finds numerous applications in calculus. Additionally, it is particularly useful in computing and representing the scattering amplitude for Regge trajectories in quantum physics, i.e., the probability amplitude of the outgoing spherical wave relative to the incoming plane wave in a stationary-state scattering process~\cite{collins,zettili,mueller}.
\begin{figure}[htbp]
\begin{center}
\begin{tikzpicture}
\draw (-4.5,0)--(0,0);
\draw[dashed,->] (0,0)--(5,0) node[right]{$x$};
\draw[->] (0,-3.5)--(0,4) node[above]{$y$};
\draw[black,fill=white] (0,0) circle (2pt);
\draw[color=blue,thick] (3,0.05) arc (0:358:3);
\draw[color=blue,thick] (0.5,0.05) arc (0:350:0.5);
\draw[color=blue,thick] (0.5,0.05)--(3,0.05);
\draw[color=blue,thick] (0.5,-0.05)--(3,-0.05);
\node[below] at (3.2,0) {$R$};
\node[below] at (-3.3,0) {$-R$};
\node[below] at (0.6,0) {$\rho$};
\draw[-{Stealth[length=1mm,width=2mm]},color=blue] (2,2.27)--(1.96,2.3);
\draw[-{Stealth[length=1mm,width=2mm]},color=blue] (-1.97,2.3)--(-2,2.27);
\draw[-{Stealth[length=1mm,width=2mm]}] (1, 0.2)--(2, 0.2);
\draw[-{Stealth[length=1mm,width=2mm]}] (2,-0.2)--(1,-0.2);
\draw[-{Stealth[length=1mm,width=2mm]},color=blue] (-0.4,0.32)--(-0.38,0.34);
\node[circle,fill,inner sep=1pt,color=red] at (-2,0){};
\node[below] at (-2.1,0) {$-1$};
\node[color=blue] at (3,1.5) {$C_R$};
\node[color=blue] at (-0.7,0.6) {$C_{\rho}$};
\end{tikzpicture}
\end{center}
\caption{The contour along a branch cut.} 		\label{branchcut}
\end{figure}
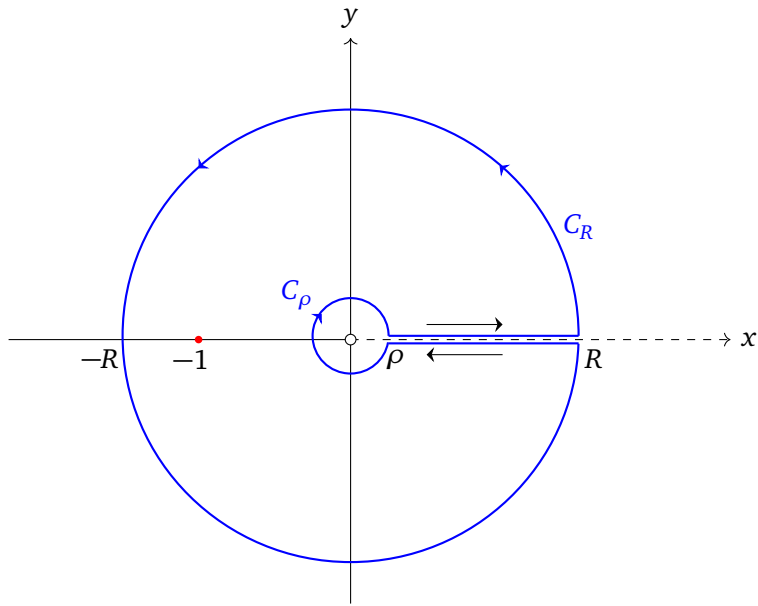

Beyond Calculus, the definite integral can be evaluated without \emph{wxMaxima} using contour integration and the Residue Theorem from Complex Analysis~\cite{boas,brown}. Let $I$ denote the definite integral, then letting $u = 1/x$, $dx = -du/u^2$ and it can be written as
\begin{equation*}
I = \int_{0}^{1} \frac{1}{(x^2 - x^3)^{\frac{1}{3}}} \, dx = \int_{\infty}^{1} \frac{-du/u^2}{(1/u^2 - 1/u^3)^{\frac{1}{3}}} = \int_{1}^{\infty} \frac{du}{u(u - 1)^{\frac{1}{3}}}.
\end{equation*}
Substitute $v = u - 1$, $dv = du$, we can write $I$ further as
\begin{equation*}
I = \int_{1}^{\infty} \frac{du}{u(u - 1)^{\frac{1}{3}}} = \int_{0}^{\infty} \frac{v^{-\frac{1}{3}}}{v + 1} \, dv.
\end{equation*}
Here, $v^{-1/3}$ denotes the positive real number or the principal value of exp$\left(-\frac{1}{3}\ln v \right)$. The final integral above is improper not only due to its upper limit of integration but also because its integrand has an infinite discontinuity at $v = 0$. Let $C_R$ and $C_{\rho}$ denote the quasi circles $|z| = R$ and $|z| = \rho$, respectively, where $\rho < 1 < R$. Figure~\ref{branchcut} displays an orientation of the contour. It is traced out by moving from $x = \rho$ to $x = R$ along the top of the branch cut for $f(z)$, then around the larger quasi circle $C_R$ and back to $x = R$, next along the bottom of the branch cut for $f$ to $x = \rho$, and finally around the smaller quasi circle $C_{\rho}$ back to $x = \rho$. We integrate the branch 
\begin{equation*}
f(z) = \frac{z^{-\frac{1}{3}}}{z + 1}, \qquad |z| > 0, \quad 0 < \text{arg} \; z < 2\pi
\end{equation*}
of the multiple-valued function $z^{-\frac{1}{3}}/(z + 1)$ around this simple closed contour with branch cut arg $z = 0$. For $z = r e^{i \theta}$, we can write
\begin{equation*}
f(z) = \frac{z^{-\frac{1}{3}}}{z + 1} = \frac{e^{-\frac{1}{3} \log z}}{z + 1} = \frac{e^{-\frac{1}{3}\left(\ln r + i \theta \right) }}{r e^{i \theta} + 1}.
\end{equation*}
On the upper edge, where $z = re^{i 0}$, we have
\begin{equation*}
f(z) = \frac{e^{-\frac{1}{3}\left(\ln r + i 0 \right) }}{r e^{i 0} + 1} = \frac{r^{-\frac{1}{3}}}{r + 1}.
\end{equation*}
On the lower edge, where $z = re^{i 2\pi}$, we obtain
\begin{equation*}
f(z) = \frac{e^{-\frac{1}{3}\left(\ln r + i 2pi \right) }}{r e^{i 2\pi} + 1} = \frac{r^{-\frac{1}{3}} e^{-i \frac{2}{3} \pi} }{r + 1}.
\end{equation*}
Applying the Residue Theorem yields
\begin{align*}
\int_{\rho}^{R} \frac{r^{-\frac{1}{3}}}{r + 1} \, dr + \int_{C_R} f(z) \, dz + \int_{R}^{\rho} \frac{r^{-\frac{1}{3}} e^{-i \frac{2}{3} \pi}}{r + 1} \, dr
+ \int_{C_{\rho}} f(z) \, dz &= 2\pi i \Res_{z = -1} f(z)		\\
\left(1 - e^{-i \frac{2}{3}\pi}\right) \int_{\rho}^{R} \frac{r^{-\frac{1}{3}}}{r + 1} \, dr &= 2\pi i e^{-i\frac{\pi}{3}} - \int_{C_R} f(z) \, dz - \int_{C_{\rho}} f(z) \, dz.
\end{align*}
We observe that
\begin{equation*}
\left|\int_{C_{\rho}} f(z) \, dz \right| \leq \int_{C_{\rho}} \left|\frac{z^{-\frac{1}{3}}}{z + 1}\right| \, dz \leq \int_{0}^{2\pi} \frac{\rho^{-\frac{1}{3}}}{1 - \rho} \, dz = \frac{2\pi}{1 - \rho} \rho^{\frac{2}{3}} \quad \to 0 \qquad \text{as} \quad \rho \to 0
\end{equation*}
and
\begin{equation*}
\left|\int_{C_{R}} f(z) \, dz \right| \leq \int_{C_{R}} \left|\frac{z^{-\frac{1}{3}}}{z + 1}\right| \, dz \leq \int_{0}^{2\pi} \frac{R^{-\frac{1}{3}}}{R - 1} \, dz = \frac{2\pi R}{R - 1} \, \frac{1}{R^{\frac{1}{3}}} \quad \to 0 \qquad \text{as} \quad R \to \infty.
\end{equation*}
By letting $\rho \to 0$ and $R \to \infty$, we arrive at
\begin{align*}
\left(1 - e^{-i \frac{2}{3}\pi}\right) \int_{0}^{\infty} \frac{r^{-\frac{1}{3}}}{r + 1} \, dr &= 2\pi i e^{-i\frac{\pi}{3}} \\
\int_{0}^{\infty} \frac{r^{-\frac{1}{3}}}{r + 1} \, dr &= 2\pi i \frac{e^{-i\frac{\pi}{3}}}{1 - e^{-i \frac{2}{3}\pi}}
= 2\pi i \frac{1}{e^{i \frac{\pi}{3} - e^{-i \frac{\pi}{3}}}} = \frac{2\pi i}{2i \sin \frac{\pi}{3}} = \frac{\pi}{\sin \frac{\pi}{3}}  = \frac{2\pi}{\sqrt{3}}.
\end{align*}
This is the same as the integral $I$ and thus we verified the computational result.

\begin{example}[Unable to solve]
Improper definite integral with infinite limits of integration
\begin{equation*}
\int_{-\infty}^{\infty} \frac{1}{(x^2 + 1)(e^x + 1)} \, dx.
\end{equation*}
While \emph{wxMaxima} provides no solution for this integral, another CAS, e.g., \emph{Mathematica}, gives the output: $\frac{\pi}{2}$.
\end{example}

\section{Methodology} \label{methodology}

We consider the qualitative aspect of the study by adopting the observational research method. In particular, we are interested in obtaining students' opinions regarding their experience after using the CAS. 

\subsection{Participant}

The participants in this study are the students who were enrolled in two Calculus courses offered at the Natural Science Campus at Sungkyunkwan University: Single Variable Calculus (SVC, course code GEDB001) and Multivariable Calculus (MVC, course code GEDB002) for three years from Spring 2016 until Fall 2019, except for the 2017 Academic Year. Each course includes three distinct cohorts, and thus we considered six sections in total. The total number of participants is 284 and their age ranges from 18 to 20 years old. There were 193 students registered in SVC and 91 were in MVC. For unclear reasons, the number of registered students for MVC tends to be fewer than for SVC. The length of the semester is 16 weeks, where Weeks~8 and~16 were designated for the Midterm and Final examination periods. We employ the convenience sampling technique in selecting the course section due to accessibility and efficiency. 

\subsection{Measurement}

We obtain the students' feedback from the emails that they sent as well as from the online questionnaire administered by the Academic Affairs Team. Students' comment are accessible through an internal network university system for faculty and students, known as \emph{Advanced Sungkyunkwan Information Square--Gold Lawn Square (ASIS-GLS)}. The students have the opportunity to deliver feedback twice: before the Midterm test and Final examination periods, i.e., Week~7 and Week~15, respectively. In particular, we pay attention to the second part of the questionnaire where the students can write freely in giving feedback to their instructors' teaching. They need to respond to the inquiry ``Please write down your suggestions for the professor to improve the class.'' The first part of the student teaching evaluation is a typical five-Likert item questionnaire and we do not cover it here since the statements are generic and do not solicit particular students' opinion on the use of CAS \emph{wxMaxima} during teaching. 

\subsection{Result} \label{result}

The response rate to the aforementioned open-ended inquiry ranges from 69\% to 98\%, depending on the semester and course. However, the total mean percentage of response rate is around 83\%. There are only 14 selected comments that are about or closely related to \emph{wxMaxima}. Table~\ref{table-participants} summarizes the details of the number of responses in each course.
\setlength{\tabcolsep}{0.5cm}
\renewcommand{\arraystretch}{1.0}
\begin{table}[h]
\caption{The distribution of participants according to the courses and number of responses.} \label{table-participants}
\begin{center}
\begin{tabular}{@{}lccrrcc@{}}
\toprule
\multirow{2}{1cm}{Course} & Number of      & \multicolumn{3}{c}{Number of responses} & \multirow{2}{1cm}{Total} & Percentage 	\\ \cline{3-5}
						  & sections 	   & Selected 			 & Midterm 			 & Final 		   &  		& average 		\\ \hline
				   SVC    & 3              & 7 				 	 & 137               & 183             & 320 	& 82.90\% 		\\
  		           MVC    & 3              & 8					 &  64               &  89             & 153 	& 84.07\%		\\ \cline{2-7}
			       Total  & 6              & 15				 	 & 201               & 272             & 473 	& 83.27\%		\\
\bottomrule       
\end{tabular}	
\end{center}
\end{table}

Regarding students' feedback, although some wrote positive comments, generally the perception of the CAS was negative. Table~\ref{feedback} displays selected students' feedback teaching evaluation related to CAS \emph{wxMaxima} for both SVC and MVC. 
\setlength{\tabcolsep}{0.25cm}
\renewcommand{\arraystretch}{1.0}
\begin{table}[h]
\caption{Selected student's feedback on teaching evaluation for SVC and MVC. We include comments related to the CAS \emph{wxMaxima} only.} \label{feedback}
\begin{center}
\begin{tabular}{@{}llll@{}}
\toprule	
Course 			& Semester    	& Period  & Feedback \\ \hline
\multirow{12}{1cm}{SVC}	& \multirow{6}{2cm}{Spring 2016} & \multirow{6}{1cm}{Midterm}		  & \emph{Maxima} program is too hard. \\ \cline{4-4}
       			&             	&         & It's nice to use \emph{Maxima}. It would be nice if there was no \\
       			&		     	&         & error by telling us how to use it or the download path. \\ \cline{4-4} 
       			&  				&  		  & The purpose of using \emph{Maxima} is good, but students may be \\
       			&             	&         & embarrassed to use the program they are not familiar with,\\
 			    &             	&         & so please provide more information on that part. \\ \cline{4-4}
 			    &				&		  &	Most of us don't know how to use those programs. \\
 			    &				&		  & With no help to learn a program, I think it is too hard and   \\
 			    &				&		  & it needs much time. \\	 \cline{3-4}
		       	&             	& Final   & \emph{Maxima} is too hard. \\ \cline{2-4}
       			& \multirow{2}{2cm}{Spring 2018}   & \multirow{2}{1cm}{Midterm} & It would be appreciated if you could reduce the drawing of \\ 
       			&				&		  &	the graph. \\ \cline{2-4}
       			& 				& 		  & Please tell me carefully the calculus.\\
       			& Spring 2019	& Midterm & It is good for taking classes and organizing concepts slowly.\\
		       	&				&		  & You show them all in a graph. \\ \hline
\multirow{14}{1cm}{MVC}&				&		  & Using \emph{Maxima}. Professor Karjanto teaches us well so I can\\
       			& 				&  		  & understand the contents well. \\ \cline{4-4}
       			&			    &         & \emph{WXMAXIMA} \\ \cline{4-4} 
       			& Fall 2016		& Midterm & \emph{maxima} program \\ \cline{4-4}
       			&				&		  & Using computer program to show graph is good.	\\ \cline{4-4}       
       			&			    &         & I want more homework solve good problem rather than \\
			    &				&		  & using computer. \\ \cline{2-4}
       			& \multirow{2}{2cm}{Fall 2018}		& \multirow{2}{1cm}{Final}	  & It seems that using more graphing software will help you \\
       			&			    &		  & understand.	\\ \cline{2-4}
       			& 				& 		  & The professor directly draws a graph to understand the \\
       			&			    & Midterm & concept and uploads additional helpful materials to  \\
	   			& Fall 2019	 	&		  & \emph{iCampus.} Was a great help in studying the subject overall. \\ \cline{3-4}
	   			&			    & \multirow{2}{1cm}{Final}   & I cannot solve the assignment using \emph{wxMaxima} and thus, \\
	   			&				&         & solve it by hand. \\
\bottomrule	
\end{tabular}
\end{center}
\end{table}

\section{Discussion} \label{discussion}

\subsection{Didactic application}

There exist some didactic applications when blending \emph{wxMaxima} with Calculus done properly. First, \emph{wxMaxima} plays a tremendous role in verifying computational results done manually by pen and paper. The opposite is also didactically useful, as \emph{wxMaxima} often spits out non-simplified outputs and this can train students to continue simplifying by hand. The former is essential in saving time, which can be used for course preparation, understanding deeper ideas, fathoming theoretical concepts, and indulging in problem-solving sessions~\cite{velychko}.

Second, obtaining the corresponding numerical values of exact expressions is as simple as inputting the command \verb|float(%)| or \verb|numer|. This can be useful when one wishes to get a sense of numerical values whereby a combination of rational numbers is an output.

Third, as expressed positively by several students, curve sketching and graphical plots generated by \emph{wxMaxima} enhance their geometrical imagination and comprehension. Appreciatively, \emph{wxMaxima} produces high-quality plots and yet lightweight for further use. Remarkably, three-dimensional plots can be rotated easily without causing any computer memory problem, as often encountered in other software, e.g., \emph{Matlab}.

Fourth, beyond Calculus, \emph{wxMaxima} has plenty of room for further exploration. It can manage data visualization and represent them with marvelous plots. It can also administer data fitting to either a straight line or nonlinear curve. It can handle programming features, such as loops, iterations, and decision making~\cite{senese}. Overall, it is useful software for those who would like to try new things.

Indeed, the literature offers abundant didactic benefits of infusing CAS in mathematics classrooms. The fifth educational application comes from a programming point of view. Coding in \emph{wxMaxima} flows logically and thus allows students to carry out simple algorithms independently~\cite{garcia}. 

Sixth, since students are allowed, even encouraged, to make mistakes and induce judgments by themselves, embracing \emph{wxMaxima} stimulates the interactive learning process through testing, evaluation, decision-making, and error correction~\cite{zakova}. Additionally, a wise application of CAS could foster students' ability in proving, modeling, problem-solving, and communication~\cite{weigand}. \emph{wxMaxima} can also be integrated with an active learning methodology such as flipped learning~\cite{karjanto19}.

\subsection{Lessons from students' feedback}

From the students' feedback, we observe that in general, the response toward CAS \emph{wxMaxima} was negative. Two comments from Spring 2016 are nearly identical, mentioned that the software is too hard to use. It is unclear whether this particular student has difficulty in downloading and installing the software package or implementing the syntax. In either case, the learning might be hampered because the CAS does not enhance the learning. For some students, the software installation was not successful, and thus they cannot operate the CAS regularly, as indicated by one comment from SVC in Spring 2016. Another student from the same cohort also commented that although using the software is recommendable, some students might be reluctant to use \emph{wxMaxima} since they are not familiar with it. 

Although we might assume that Korean digital natives grow up with technology, that does not mean that they have seen and used CAS during their previous stage of education. As one student admitted, most of them did not use any mathematical software to solve mathematical problems during their middle and high school period. Indeed, we observe a chain of logical reasoning in this situation. In order to embrace \emph{wxMaxima}, they need to attempt some problems to solve using the CAS. In order to solve Calculus problems using \emph{wxMaxima}, they need to learn how to operate and write the syntax input properly. In order to explore \emph{wxMaxima} successfully, they need proper guidance and assistance in administering the software. Thus, as instructors, we can be a guide on their side in facilitating students in discovering the usefulness of the CAS to enhance their understanding of learning Calculus.

From the selected comments, a natural conclusion would also be that more up front instruction on how to use the software and more support for using it throughout the course might improve students' dispositions toward its use. In addition to providing detailed technical operation, we could demonstrate each step from the very beginning to enact a smooth user-experience: where to download an executable file, how to install it, and how to launch the software accordingly. This is to ensure that both \emph{Maxima} and \emph{wxMaxima} are installed properly and successfully. After this essential step, we could display the most basic operations that \emph{wxMaxima} can perform, like in the same way we use a calculator to perform arithmetic problems. The tutorial might continue by showing how to evaluate a floating-point approximation, evaluating elementary Calculus problems, and sketching some simple graphs. By doing this, we hope that our digital natives will be more willing in embracing new technology.

The reaction during Fall 2016 for MVC is mixed. Some commended the use of \emph{wxMaxima} while others prefer the traditional way of solving homework by pen and paper instead of using the software. In the subsequent academic years, i.e., 2018 and 2019, what the students picked up from embedding \emph{wxMaxima} into Calculus classroom are mostly related to graph sketching. For some unclear reasons, there was no feedback related to \emph{wxMaxima} during the 2017 academic year. Thus, we excluded any feedback given in 2017 from this study. 

As the famous adage says that a picture is worth a thousand words, students can enhance their understanding of Calculus concepts by looking at figures instead of just memorizing formulas and looking at the long derivations. The feedback from the last two cohorts suggests that more students appreciate the visualization aid that \emph{wxMaxima} contributes. By sketching and displaying the graphs produced by \emph{wxMaxima}, students can comprehend better the behavior and characteristics of a particular function, whether it be 2D or 3D. The graphs help in checking whether the hand or computer calculation is correct. They also train students in observing patterns and irregularities, which in turn connect with important theorems in Calculus. By possessing this type of geometric imagination, students will be able to manipulate geometric objects in their minds, e.g., by axis rotation, scaling, shifting, etc. See~\cite{tall,mcgee,sheikh} for the role of visualization in Calculus teaching.

\subsection{Limitation} \label{limitation}

This study admits several limitations. First, we only consider the qualitative aspect of students' opinions based on an open-ended inquiry of students' feedback on teaching evaluation. We haven't investigated students' perceptions quantitatively. Admittedly, quantitative measurements might deliver better and more objective results than a handful of qualitative data. Furthermore, it would be a promising outcome to examine whether students' attitude towards \emph{wxMaxima}, or any other CAS in general, might reflect their value, self-confidence, enjoyment, motivation, and anxiety in learning Calculus or mathematics subjects. It will be equally interesting also to investigate whether the use of CAS will significantly improve students' academic performance.

Second, the course content in relation to the time availability. Prior to the 2020 academic year, we have a 16-week semester, but only 14 weeks can be used for effective teaching. Both Weeks~8 and~16 are reserved for the midterm and final examinations, respectively. For MVC, the coverage material seems to be reasonable and instructors do not need to rush in covering the material. We need to cover three chapters before the midterm and the remaining two chapters are for the second half of the semester. Thus, embedding \emph{wxMaxima} into the course seems a promising attempt. On the other hand, the materials in SVC is extremely vast while the time is insufficient. We need to encompass five chapters in the first half of the semester and the remaining four chapters need to be completed before the final exam period. Thus, dedicating one or two sessions for the \emph{wxMaxima} tutorial can be risky timewise and can jeopardize the teaching plan. Hence, splitting the course syllabus into two separate courses might offer better visibility for embedding \emph{wxMaxima} into Calculus teaching and learning.

Third, \emph{Maxima} employs command-line applications and thus can be a little bit hard to use for beginners. On the other hand, the GUI \emph{wxMaxima} can be quite helpful if one does not wish to depend on the command lines entirely. Furthermore, as we have observed in Section~\ref{symbolic}, \emph{wxMaxima} is not perfect software. When we compare with other commercial CAS like \emph{4M}, \emph{wxMaxima} is far from being ideal. A group of open-source developers who work on improving \emph{wxMaxima} is relatively smaller in numbers in comparison to professionals programmers and scientists who work in large software companies like \emph{Wolfram} or \emph{Maplesoft}. The latter work full time with handsome salaries while the former are heroes who volunteer during their spare time, sometimes as retirees. Given this kind of dramatic differences in both manpower and financial resources, we foresee that \emph{wxMaxima} will continue to fall behind other giant commercial software. As for now, the odds of any chance of catching up are still stacked against \emph{wxMaxima}.

\section{Conclusion}	\label{conclusion}

In this article, we have covered some aspects of CAS \emph{wxMaxima} where it can handle Calculus-related problems easily, with some manipulations, and unsuccessfully. Despite this imperfection, we could still utilize \emph{wxMaxima} in mathematics classrooms delicately. But, can we replace far superior software like \emph{Maple} or \emph{Mathematica} for teaching and learning? If the main concern is cost-related, then the answer to this very question would be affirmative. Yet, even if budget is not really an issue, attempting and experimenting with a new CAS might still be worth considering, at least for Calculus-linked courses, i.e., PreCalculus, SVC, and MVC. Adopting it will not only allow for both instructors and students alike to explore and learn different CAS but also to identify where \emph{wxMaxima} is particularly powerless and limited. Experiencing these situations and discovering alternative solutions can be a wonderful learning process by itself.

After an encounter with \emph{wxMaxima}, our qualitative study suggests that Korean digital natives pose a diverse reaction. Although some positively acknowledged the advantage \emph{wxMaxima} in enhancing their understanding, particularly the plotting features, some questioned whether such an attempt was necessary for Calculus teaching and learning. Indeed, several students preferred the traditional way of solving problems from textbooks with pen and paper instead of using \emph{wxMaxima}. Being freshmen at the university, the Korean digital natives might not have prior experience in using CAS when solving mathematics problems. Additionally, the lack of proper introduction might hinder the software integration process. In either case, we as instructors perceived a general reaction of amiable refusal in embracing a new technological tool. Hence, our initial assumption of digital natives readily welcoming new technology should be addressed properly and investigated further. 

After completing this study, we are longing to stimulate an accelerating debate not only among mathematics education but also the symbolic computation communities. For the former, a persuasive yet efficient way of integrating \emph{wxMaxima} into mathematics classrooms will be highly sought-after. For the latter, improving the software's effectiveness and capability in handling symbolic computations will certainly be welcomed by many researchers and educators.

\section*{Conflict of Interest}
The author declares no conflict of interest.

\section*{Dedication}
The author would like to dedicate this article to his late father Zakaria~Karjanto (Khouw~Kim~Soey, 許金瑞) who not only taught him the alphabet, numbers, and the calendar in his early childhood but also cultivated the value of hard work, diligence, discipline, perseverance, persistence, and grit. Karjanto Senior was born in Tasikmalaya, West~Java, Japanese-occupied Dutch~East~Indies on 1~January~1944 (Saturday~Pahing) and died in Bandung, West~Java, Indonesia on 18~April~2021 (Sunday~Wage).

\end{document}